\documentclass{llncs}
\usepackage{amssymb} 
\usepackage{epsfig}

\begin{document}
\pagestyle{plain}  
\bibliographystyle{splncs}

\title{Schr\"{o}dinger Equation\\As a General Optimization Algorithm}

\author{Xiaofei Huang}

\institute{eGain Communications\\Mountain View, CA 94043, U.S.A. \\
\email{huangxiaofei@ieee.org}}

\maketitle
%

\begin{abstract}
One of the greatest scientific achievements of physics in the 20th century is the discovery of quantum mechanics.
The Schr\"{o}dinger equation is the most fundamental equation in quantum mechanics describing
	the time-based evolution of the quantum state of a physical system.
It has been found that the time-independent version of the equation
	can be derived from a general optimization algorithm.
Instead of arguing for a new interpretation and possible deeper principle for quantum mechanics,
	this paper elaborates a few points of the equation as a general global optimization algorithm.
Benchmarked against randomly generated hard optimization problems, 
	this paper shows that the algorithm significantly outperformed a classic local optimization algorithm.
The former found a solution in one second with a single trial 
	better than the best one found by the latter around one hour after one hundred thousand trials.
\end{abstract}

\section{Introduction}
Optimization is a core problem both in mathematics and computer science.
It is a very active research area 
	with many international conferences every year, a large amount of literature, 
	and many researchers and practitioners across many fields for a wide range of applications. 
Combinatorial optimization~\cite{CPapadimitriou98,Pardalos02} is a branch of optimization
	where the set of feasible solutions of problems is discrete, countable, and of a finite size.
The general methods for combinatorial optimization are
   1) local search~\cite{Pardalos02,Michalewicz02}, 
   2) simulated annealing~\cite{Kirkpatrick83,Geman84}, 
   3) genetic algorithms~\cite{Fogel66,Holland92,SchwefelES}, 
   4) ant colony optimization~\cite{Dorigo2004},
   5) tabu search~\cite{Glover97}, 
   6) branch-and-bound~\cite{Lawler66,Coffman76}, 
   and  7) dynamic programming~\cite{Coffman76}.
The successful applications of different combinatorial optimization methods 
	have been reported in solving a large variety of optimization problems in practice.

Cooperative optimization is a newly proposed general global optimization method~\cite{Huang03Greece,HuangBookCCO,HuangDAGM04}
	inspired by cooperation principles in team playing.
Often times, individuals working together as a team
	can solve hard problems beyond the capability of any individual in the team.
In its normal form, cooperative optimization has one and only one equilibrium
	and converges to it with an exponential rate regardless of initial conditions.
It is also capable of identifying global optimal solutions so that it can efficiently terminate its search process.

Recently, it has been found (preprint:http://arxiv.org/abs/quant-ph/0605220) that 
	the time-independent Schr\"{o}dinger equation
	can be derived from a general form of cooperative optimization
	in a continuous-time version for continuous variables.
In particular, when the cooperation among the agents is controlled at a certain level,
	the cooperative optimization algorithm becomes the time-independent Schr\"{o}dinger equation
	when it reaches an equilibrium (a stationary state).
Instead of arguing for a new interpretation and possible deeper principle for quantum mechanics,
	this paper elaborates a few points of the equation as a general global optimization algorithm,
	such as the existence of an equilibrium and the time-based evolution.
	
In the HDTV research lab, at Tsinghua University, Beijing, China,
	we found a noticeable improvement of the algorithm over the belief propagation~\cite{Pearl88,Kschischang01}, 
	the state-of-the-art algorithm,
	at decoding some modern channel codes (a NP-hard problem similar to the integer programming problem).
At solving those real-world hard optimization problems, 
	classic optimization methods, including simulated annealing and genetic algorithms,
	failed to deliver a satisfying performance.
Instead of using the problems biased towards some particular applications,
	this paper conducts the comparison over randomly generated hard optimization problems
	to check the generality of the Schr\"{o}dinger equation as an optimization algorithm.

\section{Cooperative Optimization and Schr\"{o}dinger Equation}	

The fundamental difference between cooperative optimization and many classical optimization methods
	is at the very core of optimization, 
	i.e., the way of making decisions for assigning decision variables.
Classic one often times make a precise decision at assigning a variable 
	at every time-instance of optimization, such as $x=3$ for a time instance $t$.
Such an assignment is precise at the sense that $x$ can only be of the value $3$, not any other ones.
In contrast, cooperative optimization makes a soft decision, represented by a probability-like function called an assignment function, 
	such as $\Psi(x, t)$, at every time-instance $t$.
It says that at the time instance $t$, the variable $x$ can be of any value 
	with the likelihood measured by the function value $\Psi(x, t)$.
A variable value of a higher function value is more likely to be assigned to the variable $x$
	than anyone of a lower function value.

If the function $\Psi(x, t)$ at time $t$ is peaked at a specific value, say $x=3$, 
	while the others equal to zero,
	then the soft decision falls back to the classic precise decision, 
	e.g., assigning the value $3$ to the variable $x$ ($x=3$).
Hence, the soft decision making is a generalization of the classic precise decision making.	

Let $E(x_1, x_2, \ldots, x_n)$, simply $E(x)$, be a multivariate objective function of $n$ variables. 
Assume that $E(x)$ can be decomposed into $n$ sub-objective functions $E_i(x)$, 
	one for each variable, such that those sub-objective functions satisfying
\[ E_1(x) + E_2(x) + \ldots + E_n(x) = E(x) \ , \]
and/or the minimization of $E_i(x_i)$ with respect to $x_i$ also leads to the minimization of $E(x)$ for any $i$.

In terms of a multi-agent system, 
	let us assign $E_i(x)$ as the objective function for agent $i$, for $i=1,2,\ldots, n$.
There are $n$ agents in the system in total.
The objective of the system is to minimize $E(x)$ and 
	the objective of each agent $i$ is to minimize $E_i(x)$.

A simple form of cooperative optimization is defined as an iterative update of
	the assignment function of each agent as follows:
\begin{equation}
\Psi_i (x_i, t) = \sum_{\sim x_i} \left( e^{-E_i(x)/\hbar} \prod_{j \not= i} p_j(x_j, t-1) \right),\quad \mbox{for $i=1,2,\ldots,n$}  \ , 
\label{cooperative_optimization_general3}
\end{equation}
where $\sum_{\sim x_i}$ stands for the summation over all variables except $x_i$ and $\hbar$ is a constant of a small positive value.
$p_i(x_i, t)$ is defined as 
\begin{equation} 
p_i(x_i, t) = \left(\Psi_i(x_i,t)\right)^{\alpha} / \sum_{x_i} \left(\Psi_i(x_i,t)\right)^{\alpha} \ , 
\label{compute_assignment_probabilty}
\end{equation}
where $\alpha$ is a parameter of a positive real value.

By the definition, $p_i(x_i, t)$ just likes a probability function satisfying
\[ \sum_{x_i} p_i (x_i, t) = 1 \ . \]
It is, therefore, called the assignment probability function.
It defines the probability-like soft decision at assigning variable $x_i$ at the time instance $t$.

The original assignment function $\Psi_i(x_i, t)$, is called the assignment state function.
That is, the state of agent $i$ at the time instance $t$ is represented 
	by its assignment state function $\Psi_i(x_i, t)$.
From Eq.~\ref{compute_assignment_probabilty} we can see that
	the assignment probability function $p_i(x_i, t)$ is defined 
	as the assignment state function $\Psi_i(x_i)$ to the power $\alpha$ with normalization.
To show the relationship, the assignment probability function $p_i(x_i, t)$ is also expressed as $\left({\bar \Psi}_i(x_i,t)\right)^{\alpha}$ 
	in the following discussions with the bar standing for the normalization.
	
With this notation, the iterative update function~(\ref{cooperative_optimization_general3}) can be rewritten as
\begin{equation}
\Psi_i (x_i, t) = \sum_{\sim x_i} \left( e^{-E_i(x)/\hbar} \prod_{j \not= i} \left({\bar \Psi}_j (x_j, t-1)\right)^{\alpha} \right),\quad \mbox{for $i=1,2,\ldots,n$}  \ . 
\label{cooperative_optimization_general3b}
\end{equation}

Without loss of generality, let $e^{-E_i(x)/\hbar}$ be the utility function of agent $i$ in terms of game theory.
It is important to note that when the parameter $\alpha$ is of a sufficient value, i.e.,
	$\alpha \rightarrow \infty$,
	the simple form~(\ref{cooperative_optimization_general3b})
	converges to an equilibrium if and only if it is also a Nash equilibrium (preprint: http://arxiv.org/abs/0901.3615).
Nash Equilibrium is arguably the most important concept in game theory,
	critical for understanding a common scenario in game playing.
It offers the mathematical foundation for social science and economy.

By substituting Eq.~\ref{cooperative_optimization_general3} into Eq.~\ref{compute_assignment_probabilty},
	we have a mapping from a set of assignment probability functions to itself.
Because the set is compact and the mapping is continuous, 
	so a fixed point exists based on Brouwer fixed point theorem.
Since a set of assignment state functions 
	is uniquely defined by a set of assignment probability functions by Eq.~\ref{cooperative_optimization_general3},
We can conclude that 	
	there exists at least one set of assignment state functions $\{\Psi^{*}_1(x_1),\Psi^{*}_2(x_2), \ldots, \Psi^{*}_n(x_n) \}$ such that
\[\Psi^{*}_i (x_i) = \sum_{\sim x_i} \left( e^{-E_i(x)/\hbar} \prod_{j \not= i} \left({\bar \Psi}^{*}_j(x_j)\right)^{\alpha} \right),~~~\mbox{for $i=1,2,\ldots, n$} \ . \]
	
In particular, when $\alpha = 2$, 
	the simple general form of cooperative optimization over the real domain $\mathbb{R}$ given in (\ref{cooperative_optimization_general3b})
	can be generalized further over the complex domain $\mathbb{C}$ in a continuous-time version as follows
\begin{equation}
\hbar \frac{\partial \psi_i (x_i, t)}{\partial t} = -\frac{1}{Z_i (t)} \psi_i (x_i, t) \sum_{\sim x_i} \left( E_i (x) \prod_{j \not= i} |\psi_j(x_j, t)|^2  \right) \ , 
\label{cooperative_optimization_general3.2}
\end{equation}
where $Z_i(t)$ is a normalization factor such that 
\[ \sum_{x_i} |\psi_i(x_i, t)|^2 = 1 \ , \]
and $|\psi_i(x_i, t)|^2$ is defined as
\[ |\psi_i(x_i, t)|^2 = \psi^{*}_i(x_i, t) \psi_i(x_i, t) \ . \]

Following the notation from physics, denote $\psi_i (x_i, t)$ and a vector $ \mid \psi_i (t) \rangle$.
The equation~(\ref{cooperative_optimization_general3.2}) can be generalized further as follows
\begin{equation}
\hbar \frac{d}{d t} \mid \psi_i (t) \rangle = -\frac{1}{Z_i (t)} H_i \mid \psi_i (t) \rangle \ . 
\label{cooperative_optimization_general3.3}
\end{equation}
In the equation, $H_i$ is a hermitian matrix which defines the local energy (objective) function for agent $i$ in a more general form.
In particular, for the case of (\ref{cooperative_optimization_general3.2}),
	the hermitian matrix $H_i$ reduces to a diagonal matrix with the diagonal elements as 
\[ \sum_{X_i \setminus{x_i}} \left(E_i (x) \prod_{j \not= i} |\psi_j(x_j, t)|^{2} \right), \quad \mbox{for every value of $x_i$} \ . \]

The function $\psi_i (x_i, t)$ is also called a wavefunction in physics.
It is important to note that the equation~(\ref{cooperative_optimization_general3.3}) 
	is the dual equation of the Schr\"{o}dinger equation:
\[ i \hbar \frac{d}{d t} \mid \psi_i (t) \rangle = H_i \mid \psi_i (t) \rangle \ , \] 
where $-1$ is replaced by the imaginary unit $i$($=\sqrt{-1})$ and the normalization factor $Z_i(t)$ is not required
	since the equation is unitary, 
	which means that the total norm of the wavefunction $\psi_i(x_i, t)$ is reserved, i.e.,
\[ \sum_{x_i} |\psi_i(x_i, t)|^2 = const, \quad \mbox{for any $t$.} \]

When the dynamic equation~(\ref{cooperative_optimization_general3.3}) 
	reachs a stationary point (equilibrium),
	it has been found out that the equation becomes the time-independent Schr\"{o}dinger equation, one of the most important equations in quantum mechanics.
That is
\[ e_i \mid \psi_i (x_i, t) \rangle = H_i \mid \psi_i (x_i, t) \rangle \ , \]
where $e_i$ can only be any one of the eigenvalues of $H_i$.

If the local energy function $H_i$ is time-independent, 
	then the time-based evolution of the state function $\mid \psi_i (t) \rangle$ is
\[ \mid \psi_i (t) \rangle= \frac{1}{Z_i (t)} e^{-\frac{1}{\hbar} H_i t} \mid \psi_i (0) \rangle \ . \]	

Assume that $\mid \phi_1 \rangle, \mid \phi_2 \rangle, \ldots, \mid \phi_N \rangle$  are 
	the $N$ eigenvectors of the Hamilitian matrix $H_i$
	with corresponding eigenvalues as $e_1, e_2, \ldots, e_n$, respectively.
Then any initial state $\mid \psi_i (0) \rangle$ can be represented as the linear superposition
	of the $N$ eigenvectors as follows
\[\mid \psi_i (0) \rangle = c_1 \mid \phi_1 \rangle +  c_2 \mid \phi_2 \rangle + \cdots + c_N \mid \phi_N \rangle \ , \]
where $c_n$ ($1 \le n \le N$ and $c_n \in \mathbb{C}$) are coefficients.
The state function $\mid \psi_i (x_i, t) \rangle$ at time $t$ is governed by 
\[ \mid \psi_i (t) \rangle = \frac{1}{Z_i (t)} \sum^{N}_{j=1} c_j e^{-\frac{1}{\hbar}e_j t} \mid\phi_j\rangle \ . \]

Benefit from almost a century's mathematical development driven quantum physics, 
	we can extend the concept of a hermitian matrix to that of a hermitian operator 
	to deal with the case when the variable $x_i$ is continuous.
In this case, the concept of an eigenvector will be extended to that of an eigenfunction.

Putting everything together,
	the pseudo-code of a quantum optimization algorithm is shown in Figure~\ref{quantum_optimization_algorithm}.

\begin{figure}
\begin{tabbing}
123\=456\=789\=012\=345\=678\=901\ \kill \\
1\>{\bf begin}\\
2\>\> for every $i$, initialize $\Psi_i(x_i)$ with random non-negative values and normalize it; \\
3\>\>{\bf for}  $k:=1$ {\bf to} $maximum\_iteration\_number$ {\bf do} \\
4\>\> {\bf begin}\\
5\>\>\>{\bf for}  $i:=1$ {\bf to} $n$ {\bf do} \\
6\>\>\> {\bf begin}\\
7\>\>\>\>{\bf for} each value of $x_i$ {\bf do} \\
8\>\>\>\>\> $\Psi_i (x_i) := \sum_{\sim x_i} \left( e^{-E_i(x)/\hbar} \prod_{j \not= i} {\bar \Psi}^2_j(x_j) \right)$;\\
9\>\>\>\>Let $\tilde{x}_i := \arg \max_{x_i} \Psi_i(x_i)$;\\
10\>\>\>{\bf end};\\
11\>\>\>normalize $\Psi_i(x_i)$ such that $\sum_{x_i} |\psi_i(x_i)|^2 = 1$; \\
12\>\> {\bf end};\\
13\>\> return $(\tilde{x}_1, \tilde{x}_2,\ldots, \tilde{x}_n)$ as the final solution;\\
14\>{\bf end};\\
\end{tabbing}
\caption{A quantum optimization algorithm}
\label{quantum_optimization_algorithm}
\end{figure}

\section{Experimental Comparisons}	
\subsection{Constraint Optimization Problems}

A large class of optimization problems has an objective function of the following form,
\begin{equation}
E(x)=\sum^{n}_{i=1} f_i (x_{i}) + \sum_{(i,j) \in {\cal N}} f_{ij} (x_i, x_j) \ . 
\label{binary_cost_function}
\end{equation}
The function $f_i(x_i)$ is called an unary function
	and the function $f_{ij}(x_i, x_j)$ is called a binary function.
To note the collection of all defined binary functions,
	the set ${\cal N}$ is used which contains non-ordered pairs of variable indices
	where each pair $(i,j)$ corresponds to a defined binary function $f_{ij}(x_i, x_j)$.

Minimizing the objective function defined in (\ref{binary_cost_function}) is referred to as a constraint optimization problem (binary COP) in AI.
The unary function $f_i(x)$ is called an unary constraint on variable $x_i$
	and the binary function $f_{ij}(x_i, x_j)$ is called a binary constraint on variables $x_i, x_j$.

The constraint optimization problem is a very general formulation 
	for many optimization problems arose from widely different fields.
Examples are the famous traveling salesman problem,
	the weighted maximum satisfiability problem, 
	the quadratic variable assignment problem, 
	stereo matching in computer vision, 
	image segmentation in image processing, and many more.
The constraint optimization problem belongs to the class of NP-complete problems.

An objective function in form of (\ref{binary_cost_function})
	can be represented with an undirected graph.
In the graph, each variable $x_i$ is represented by a node and 
	each binary constraint $f_{i,j}(x_i, x_j)$ is represented by an undirected edge
	connecting variable node $x_i$ to variable node $x_j$.
In this graphic representation,
	the degree of a node is measured by the total number of edges connected to it.

Each variable node $x_i$ can have a number of neighboring variable nodes.
Let ${\cal N}(i)$ be the set of the indices of the neighboring variables of $x_i$. 
By definition, 
\[ {\cal N}(i) = \{j|(i,j) \in {\cal N} \} \ . \]

Using the notations, we can define a local objective function for agent $i$ as 
\begin{equation}
E_i(x)= f_i (x_{i}) + \sum_{j \in {\cal N}(i)} f_{ij} (x_i, x_j), \quad \mbox{for $i=1,2, \ldots, n$} \ .
\label{binary_cost_function2}
\end{equation}
Obviously, minimizing the local objective function $E_i(x)$ with respect to $x_i$ 
	also leads to the minimization of the global objective function $E(x)$.

\subsection{A Classic Local Optimization Algorithm}

Local search is one of the classic methods for attacking hard combinatorial optimization problems.
It also plays a fundamental role at understanding many other optimization methods,
	such as simulated annealing,
   	genetic algorithms, 
	and tabu search.
Started from a randomly generated initial solution,
	local search iteratively tries to replace the current solution 
	by a better one in the neighborhood of the current solution until no further improvement is possible.
	
The improvement of the local search algorithm choosing in the experiments is 
	achieved by adjusting the value of each variable while other variables have their values fixed.
Specifically, assume that the current solution is $\tilde{x}$ and $x_i$ is the variable with its value to be adjusted.
Its value in the current solution is $\tilde{x}_i$.
Let $x_{-i}$ denote all variables except $x_i$.
In this case, the neighborhood is defined as $(x_i, \tilde{x}_{-i})$ where $x_i \in D_i$, the domain of $x_i$.
The adjustment is achieved by finding the best value for $x_i$, denote as $x^{*}_i$, such that
\[ E(x^{*}_i) = \min_{x_i} E(x_i, \tilde{x}_{-i}) \ , \]
If $x^{*}_i \not=\tilde{x}_i$, then replacing $\tilde{x}_i$ by $x^{*}_i$.
Otherwise, keeping $\tilde{x}_i$ as the value of $x_i$.
If $x^{*}_i = \tilde{x}_i$ and it is true for every variable $x_i$, 
	then no solution improvement is possible by adjusting any variable. 
In this case, the local search algorithm is called reaching a local optimal solution 
	since it is optimal with respect only to its neighbors.
It terminates the search process and returns the current solution $\tilde{x}$ as the final one.

In particular, if $E(x)$ can be decomposed into $n$ sub-objective functions $E_i(x)$ 
	and minimizing $E_i(x)$ with respect to $x_i$ also leads to minimizing $E(x)$,
	then the above optimization problem becomes
\begin{equation}
E_i(x^{*}_i) = \min_{x_i} E_i(x_i, \tilde{x}_{-i})  \ . 
\label{local_optimization}
\end{equation}

A pseudo-code of the above local search algorithm is given in Figure~\ref{local_search}.

\begin{figure}
\begin{tabbing}
123\=456\=789\=012\=345\=678\=901\=222\ \kill \\
1\>{\bf begin}\\
5\>\> randomly generate a solution $\tilde{x}$; \\
6\>\> {\bf repeat} \\
7\>\>\> $local\_minimum\_flag :=true;$ \\
8\>\>\>{\bf for}  $i:=1$ {\bf to} $n$ {\bf do} \\
9\>\>\> {\bf begin}\\
10\>\>\>\> find $x^{*}_i$ such that $E(x^{*}_i) = \min_{x_i} E(x_i, \tilde{x}_{-i})$; \\
11\>\>\>\> {\bf if} $x^{*}_i \not = \tilde{x}_i$ {\bf then} \\
12\>\>\>\> {\bf begin}\\
13\>\>\>\>\> $\tilde{x}_i:=x^{*}_i$; \\
14\>\>\>\>\> $local\_minimum\_flag :=false$; \\
15\>\>\>\> {\bf end};\\
16\>\>\> {\bf end};\\
17\>\> {\bf until} $local\_minimum\_flag = true$;\\
20\>\>{\bf return} $\tilde{x}$ as the final solution;\\
21\>{\bf end};\\
\end{tabbing}
\caption{The local search algorithm}
\label{local_search}
\end{figure}


The multi-restart local search algorithm simply calls the above local search algorithm
	multiple times and picks the best solution (the best local optimal solution) as the final solution.
In other words, it discovers multiple local optimal solutions and keeps the best one.
	

\subsection{The Experimental Results}

The size of a randomly generated instance is controlled by the total number of variables 
	and the total number of values for each variable.
The edges connecting to each variable node are randomly selected.
If each node is connected to all other nodes, then the graph is a full graph.
Otherwise, if it is connected to a small percentage of other nodes, then the graph is sparsely connected.
Therefore, we can use the average node degree to control the graphical structure of each generated instance.

The function values of the constraints $f_i(x_i)$ and $f_{ij}(x_i, x_j)$
	are uniformly sampled from the interval $[0, 1]$.
Because of the randomness, 
	it is hard to apply any domain-specific heuristics for the advantage of any optimization algorithm.

In the first set of experiments, 
	ten instances of the constraint optimization problem are generated.
Each instance has $121$ variables, $50$ values for each variable, and the average node degree is six.
Each instance has an enormous number of feasible solutions,
	bigger than the total number of atoms in the universe,
	posing a challenging optimization problem.

For each of those ten randomly generated instances, 
	the best solution of the multi-restart local search algorithm (MRLS) after 100,000 trials
	is compared with the solution found by 
	the quantum optimization algorithm (QOA) (Fig.~\ref{quantum_optimization_algorithm}) with a single trial.
The parameters of the algorithm are set as $\hbar=1$ and $maximum\_iteration\_number =20$.
The results are given in the following table.
From it we can see that the quantum optimization algorithm yields better results within a much shorter time
	than the multi-restart local search algorithm.
In the table, $variable\#$ represents the total number of variables
	and $value\#$ the total number of values for each variable.
	
\begin{tabular}{cccc}
\hline
\multicolumn{4}{c}{MRLS with 100,000 trials {\it vs} QOA with a single trial} \\
\multicolumn{4}{c}{variable\#=121, value\#=50, average\_node\_degree=6} \\
\hline
Instance\# & MRLS:cost (time) &  QOA:cost (time) & improvement\\
\hline
1 & 153.11 (3001 sec) & 144.77 (1.061 sec) & 5.76\%\\
2 & 153.92 (3067 sec) & 144.02 (1.046 sec) & 6.87\%\\
3 & 152.12 (2994 sec) & 136.84 (1.030 sec) & 11.17\% \\
4 & 153.53 (3072 sec) & 144.34 (1.045 sec) & 6.37\% \\
5 & 154.95 (3013 sec) & 145.22 (1.061 sec) & 6.70\% \\
6 & 148.17 (3012 sec) &140.37 (1.029 sec) & 5.55\% \\
7 & 147.90 (2905 sec) &138.83 (1.014 sec) & 6.54\% \\
8 & 171.03 (3275 sec) &158.82 (1.138 sec) & 7.69\% \\
9 & 154.70 (3013 sec) &145.59 (1.061 sec) & 6.26\% \\
10 & 145.14 (2972 sec) & 130.19 (1.030 sec) & 11.49\% \\
\hline
\end{tabular}

The following table lists the experimental results with another set of instances 
	where their sizes are different from the previous ones.
The total number of variables is inceased to $1001$ and the total number of values for each variable is reduced to $10$.
From the table we can see that the quantum optimization algorithm still outperforms
	the multi-restart local search algorithm both in quality and speed.

\begin{tabular}{cccc}
\hline
\multicolumn{4}{c}{MRLS with 10,000 trials {\it vs} QOA with a single trial} \\
\multicolumn{4}{c}{variable\#=1001,value\#=10,average\_node\_degree=10} \\
\hline
Instance\# & MRLS:cost (time) &  QOA:cost (time) & improvement\\
\hline
  1 & 3269.97 (41.6 sec) & 3102.63 (1.373 sec) & 5.39\% \\
  2 & 3221.61 (4126 sec) & 3084.81 (1.415 sec) & 4.43\% \\
  3 & 3237.84 (4179 sec) & 3090.48 (1.389 sec) & 4.77\% \\
  4 & 3270.37 (4134 sec) & 3159.82 (1.435 sec) & 3.50\% \\
  5 & 3267.66 (4113 sec) & 3109.14 (1.373 sec) & 5.10\% \\
  6 & 3307.75 (4108 sec) & 3204.13 (1.404 sec) & 3.23\% \\
  7 & 3248.23 (4073 sec) & 3153.07 (1.389 sec) & 3.02\% \\
  8 & 3273.33 (4077 sec) & 3146.69 (1.388 sec) & 4.02\% \\
  9 & 3300.05 (4117 sec) & 3188.34 (1.388 sec) & 3.50\% \\
 10 & 3269.44 (4126 sec) & 3141.70 (1.454 sec) & 4.07\% \\
\hline
\end{tabular}

\section{Conclusions}
The Schr\"{o}dinger equation can be converted into a general purpose optimization algorithm.
The original equation is in a continuous-time version 
	while the one presented in this paper is in a discrete-time version
	so that it can be easily implemented by computer software or hardware.
Benchmarked against randomly generated hard optimization problems, 
	the quantum optimization algorithm significantly outperformed a classic multi-restart local search algorithm
	by several orders of magnitude.
Together with many other existing optimization methods,
	the quantum optimization can be served as a new method to attack hard optimization problems.

\nocite{Pardalos02}
\nocite{GameTheoryLuce}

\nocite{CPapadimitriou98}
\nocite{Messiah99,Tegmark01,Seife05}
\nocite{GameTheoryLuce,Daskalakis05three-playergames,DaskalakisGoldbergPapadimitriou05}


\begin{thebibliography}{10}

\bibitem{CPapadimitriou98}
Papadimitriou, C.H., Steiglitz, K., eds.:
\newblock Combinatorial Optimization.
\newblock Dover Publications, Inc. (1998)

\bibitem{Pardalos02}
Pardalos, P., Resende, M.:
\newblock Handbook of Applied Optimization.
\newblock Oxford University Press, Inc. (2002)

\bibitem{Michalewicz02}
Michalewicz, Z., Fogel, D.:
\newblock How to Solve It: Modern Heuristics.
\newblock Springer-Verlag, New York (2002)

\bibitem{Kirkpatrick83}
Kirkpatrick, Gelatt, C., Vecchi, M.:
\newblock Optimization by simulated annealing.
\newblock Science \textbf{220} (1983)  671--680

\bibitem{Geman84}
Geman, S., Geman, D.:
\newblock Stochastic relaxation, gibbs distributions, and the bayesian
  restoration of images.
\newblock IEEE Transactions on Pattern Analysis and Machine Intelligence
  \textbf{PAMI-6} (1984)  721--741

\bibitem{Fogel66}
Fogel, L., Owens, A., Walsh, M.:
\newblock Artificial Intelligence through Simulated Evolution.
\newblock John Wiley, New York (1966)

\bibitem{Holland92}
Hinton, G., Sejnowski, T., Ackley, D.:
\newblock Genetic algorithms.
\newblock Cognitive Science (1992)  66--72

\bibitem{SchwefelES}
Schwefel, H.P.:
\newblock Evolution and Optimum Seeking.
\newblock John Wiley and Sons, Inc. (1995)

\bibitem{Dorigo2004}
Dorigo, M., St\"{u}tzle, T.:
\newblock Ant Colony Optimization.
\newblock The MIT Press, Cambridge, Massachusetts, London, England (2004)

\bibitem{Glover97}
Golver, F., Laguna, M.:
\newblock Tabu Search.
\newblock Kluwer Academic Publishers (1997)

\bibitem{Lawler66}
Lawler, E.L., Wood, D.E.:
\newblock Brand-and-bound methods: A survey.
\newblock OR \textbf{14} (1966)  699--719

\bibitem{Coffman76}
Jr., E.G.C., ed.:
\newblock Computer and Job-Shop Scheduling.
\newblock Wiley-Interscience, New York (1976)

\bibitem{Huang03Greece}
Huang, X.:
\newblock A general framework for constructing cooperative global optimization
  algorithms.
\newblock In: Frontiers in Global Optimization. Nonconvex Optimization and Its
  Applications.
\newblock Kluwer Academic Publishers (2004)  179--221

\bibitem{HuangBookCCO}
Huang, X.:
\newblock Cooperative optimization for solving large scale combinatorial
  problems.
\newblock In: Theory and Algorithms for Cooperative Systems. Series on
  Computers and Operations Research.
\newblock World Scientific (2004)  117--156

\bibitem{HuangDAGM04}
Huang, X.:
\newblock Cooperative optimization for energy minimization in computer vision:
  A case study of stereo matching.
\newblock In: Pattern Recognition, 26th DAGM Symposium, Springer-Verlag, LNCS
  3175 (2004)  302--309

\bibitem{Pearl88}
Pearl, J.:
\newblock Probabilistic Reasoning in Intelligent Systems: Networks of Plausible
  Inference.
\newblock Morgan Kaufmann (1988)

\bibitem{Kschischang01}
Kschischang, F.R., Frey, B.J., andrea Loeliger, H.:
\newblock Factor graphs and the sum-product algorithm.
\newblock IEEE Transactions on Information Theory \textbf{47} (2001)  498--519

\bibitem{GameTheoryLuce}
Luce, R.D., Raiffa, H.:
\newblock Games and Decisions: Introduction and Critical Survey.
\newblock Dover (1985)

\bibitem{Messiah99}
Messiah, A.:
\newblock Quantum Mechanics.
\newblock Dover Publications, Inc., Mineola, New York (1999)

\bibitem{Tegmark01}
Tegmark, M., Wheeler, J.A.:
\newblock 100 years of quantum mysteries.
\newblock Scientific American (2001)  68--75

\bibitem{Seife05}
Seife, C.:
\newblock Do deeper principles underlie quantum uncertainty and nonlocality?
\newblock Science \textbf{309} (2005) ~98

\bibitem{Daskalakis05three-playergames}
Daskalakis, C., Papadimitriou, C.H.:
\newblock Three-player games are hard.
\newblock In: Electronic Colloquium on Computational Complexity (ECCC. (2005)
  TR05--139

\bibitem{DaskalakisGoldbergPapadimitriou05}
Daskalakis, C., Goldberg, P., Papadimitriou, C.:
\newblock The complexity of computing a nash equilibrium.
\newblock In: Electronic Colloquium on Computational Complexity (ECCC. (2005)
  TR05--115

\end{thebibliography}

\end{document}